\documentclass{ws-p8-50x6-00}

\begin{document}

\title{Inversion formulas involving orthogonal polynomials and
some of their applications}

\author{Roelof Koekoek}

\address{Delft University of Technology, Faculty of Information
Technology and Systems, P.O. Box 5031, 2600 GA Delft, The
Netherlands\\E-mail: koekoek@twi.tudelft.nl}

\maketitle

\abstracts{
We derive inversion formulas involving orthogonal polynomials
which can be used to find coefficients of differential equations
satisfied by certain generalizations of the classical orthogonal
polynomials. As an example we consider special symmetric
generalizations of the Hermite polynomials.}

\section{Introduction}
In \cite{Dvlag} we found differential equations of spectral type
satisfied by the generalized Laguerre polynomials
$\left\{L_n^{\alpha,M}(x)\right\}_{n=0}^{\infty}$ which are
orthogonal on the interval $[0,\infty)$ with respect to the
weight function
$$\frac{1}{\Gamma(\alpha+1)}x^{\alpha}e^{-x}+M\delta(x),\;\alpha>-1,\;M\ge 0.$$
These orthogonal polynomials were introduced by T.H.~Koornwinder
in \cite{Koorn}. In order to find the coefficients of these
differential equations we had to solve systems of equations of the
form
\begin{equation}
\label{sysLag}
\sum_{i=1}^{\infty}a_i(x)D^iL_n^{(\alpha)}(x)=F_n(x),\;n=1,2,3,\ldots,
\end{equation}
where $\displaystyle D=\frac{d}{dx}$ denotes the differentiation
operator. In \cite{Bavinck} H.~Bavinck showed that the
coefficients $\{a_i(x)\}_{i=1}^{\infty}$ are uniquely determined
and can be written in the form
\begin{equation}
\label{oplLag}
a_i(x)=(-1)^i\sum_{j=1}^iL_{i-j}^{(-\alpha-i-1)}(-x)F_j(x),\;i=1,2,3\ldots.
\end{equation}
This result is based on the inversion formula
\begin{equation}
\label{invLag}
\sum_{k=j}^iL_{i-k}^{(-\alpha-i-1)}(-x)L_{k-j}^{(\alpha+j)}(x)=
\delta_{ij},\;j\le i,\;i,j=0,1,2,\ldots.
\end{equation}
See also \cite{Invjac}. This inversion formula was derived in a
similar way as the inversion formula involving Charlier
polynomials found in \cite{Charlier}. See also \cite{Invjac} and
section~3 of this paper. For inversion formulas involving Meixner
polynomials the reader is referred to \cite{Meixner}. See also
\cite{Meixner2}. In \cite{Soblag} we used the inversion formula
(\ref{invLag}) to find differential equations of spectral type
satisfied by the Sobolev-type Laguerre polynomials
$\left\{L_n^{\alpha,M,N}(x)\right\}_{n=0}^{\infty}$ which are
othogonal with respect to the Sobolev-type inner product
$$<f,g>\,=\frac{1}{\Gamma(\alpha+1)}\int_{0}^{\infty}x^{\alpha}e^{-x}f(x)g(x)dx
+Mf(0)g(0)+Nf'(0)g'(0),$$
where $\alpha>-1$, $M\ge 0$ and $N\ge 0$.

In \cite{Koorn} T.H.~Koornwinder also introduced the generalized Jacobi
polynomials $\left\{P_n^{\alpha,\beta,M,N}(x)\right\}_{n=0}^{\infty}$
which are orthogonal on the interval $[-1,1]$ with respect to the
weight function
$$\frac{\Gamma(\alpha+\beta+2)}{2^{\alpha+\beta+1}\Gamma(\alpha+1)\Gamma(\beta+1)}
(1-x)^{\alpha}(1+x)^{\beta}+M\delta(x+1)+N\delta(x-1),$$
where $\alpha>-1$, $\beta>-1$, $M\ge 0$ and $N\ge 0$. In
\cite{Search} we were looking for differential equations of
spectral type satisfied by these generalized Jacobi polynomials.
The general case turned out to be very difficult, but in \cite{Symjac}
we were able to solve this problem in the special case that
$\beta=\alpha$ and $N=M$. In order to find the coefficients of
these differential equations we had to solve systems of equations
of the form
$$\sum_{i=1}^{\infty}c_i(x)D^iP_n^{(\alpha,\beta)}(x)=F_n(x),
\;n=1,2,3,\ldots.$$
In \cite{Invjac} we showed that the coefficients $\{c_i(x)\}_{i=1}^{\infty}$
are unique and that they can be written in the form
$$c_i(x)=2^i\sum_{j=1}^i\frac{\alpha+\beta+2j+1}{(\alpha+\beta+j+1)_{i+1}}
P_{i-j}^{(-\alpha-i-1,-\beta-i-1)}(x)F_j(x),\;i=1,2,3\ldots.$$
This result is based on the inversion formula
\begin{eqnarray}
\label{invJac}
& &\sum_{k=j}^i\frac{\alpha+\beta+2k+1}{(\alpha+\beta+k+j+1)_{i-j+1}}\nonumber\\
& &{}\times P_{i-k}^{(-\alpha-i-1,-\beta-i-1)}(x)
P_{k-j}^{(\alpha+j,\beta+j)}(x)=\delta_{ij},\;j\le i,\;i,j=0,1,2,\ldots,
\end{eqnarray}
which is proved in \cite{Invjac}. This inversion formula was
derived in a completely different way than the inversion formulas
mentioned before. In \cite{Madrid} it is shown that this inversion
formula (with $\beta=\alpha$) can be used to derive the results
obtained in \cite{Symjac} in an easier way. Finally in
\cite{Dvjac} this inversion formula is used to solve the problem
for all $\alpha>-1$, $\beta>-1$, $M\ge 0$ and $N\ge 0$.

In this paper we will derive several kinds of inversion formulas
and we will show how they can be applied to find coefficients of
differential equations for generalizations of some classical
orthogonal polynomials.

\section{Some classical orthogonal polynomials}

In this section we will recall some formulas involving classical
orthogonal polynomials which we will use in this paper. For
details the reader is referred to \cite{Askey}.

The Meixner-Pollaczek polynomials $\left\{P_n^{(\lambda)}(x;\phi)\right\}_{n=0}^{\infty}$
can be defined by their generating function
\begin{equation}
\label{genMP}
\left(1-e^{i\phi}t\right)^{-\lambda+ix}\left(1-e^{-i\phi}t\right)^{-\lambda-ix}
=\sum_{n=0}^{\infty}P_n^{(\lambda)}(x;\phi)t^n.
\end{equation}

The classical Jacobi polynomials $\left\{P_n^{(\alpha,\beta)}(x)\right\}_{n=0}^{\infty}$
can be defined for all $\alpha$ and $\beta$ and $n\in\{0,1,2,\ldots\}$ by
\begin{equation}
\label{defJac}
P_n^{(\alpha,\beta)}(x)=\sum_{k=0}^n\frac{(n+\alpha+\beta+1)_k}{k!}
\frac{(\alpha+k+1)_{n-k}}{(n-k)!}\left(\frac{x-1}{2}\right)^k.
\end{equation}
They satisfy the orthogonality relation
\begin{eqnarray*}
& &\frac{\Gamma(\alpha+\beta+2)}{2^{\alpha+\beta+1}\Gamma(\alpha+1)\Gamma(\beta+1)}
\int_{-1}^1(1-x)^{\alpha}(1+x)^{\beta}P_m^{(\alpha,\beta)}(x)P_n^{(\alpha,\beta)}(x)dx\\
& &{}\hspace{1cm}=\frac{\alpha+\beta+1}{2n+\alpha+\beta+1}\,
\frac{(\alpha+1)_n(\beta+1)_n}{(\alpha+\beta+1)_n\,n!}\,\delta_{mn},\;m,n=0,1,2,\ldots.
\end{eqnarray*}
The Gegenbauer or ultraspherical polynomials $\left\{G_n^{(\lambda)}(x)\right\}_{n=0}^{\infty}$
form a special case of the classical Jacobi polynomials. In fact
we have
\begin{equation}
\label{rel1}
G_n^{(\lambda)}(x)=\frac{(2\lambda)_n}{(\lambda+\frac{1}{2})_n}
P_n^{(\lambda-\frac{1}{2},\lambda-\frac{1}{2})}(x),
\;\lambda>-\frac{1}{2},\;\lambda\ne 0.
\end{equation}
These ultraspherical polynomials can also be defined by their
generating function
\begin{equation}
\label{genultra}
(1-2xt+t^2)^{-\lambda}=\sum_{n=0}^{\infty}G_n^{(\lambda)}(x)t^n.
\end{equation}
The special case $\lambda=0$ needs another normalization. In that
case we have the Chebyshev polynomials of the first kind $\left\{T_n(x)\right\}_{n=0}^{\infty}$
given by
$$T_n(x)=\frac{P_n^{(-\frac{1}{2},-\frac{1}{2})}(x)}{P_n^{(-\frac{1}{2},-\frac{1}{2})}(1)}=
{}_2F_1\left(\left.{{-n, n} \atop \frac{1}{2}}\right|\frac{1-x}{2}\right),\;n=0,1,2,\ldots.$$
Their generating function equals
\begin{equation}
\label{genCheb1}
\frac{1-xt}{1-2xt+t^2}=\sum_{n=0}^{\infty}T_n(x)t^n.
\end{equation}
The Chebyshev polynomials of the second kind $\left\{U_n(x)\right\}_{n=0}^{\infty}$
are given by
$$U_n(x)=(n+1)\,\frac{P_n^{(\frac{1}{2},\frac{1}{2})}(x)}{P_n^{(\frac{1}{2},\frac{1}{2})}(1)}=
(n+1)\,{}_2F_1\left(\left.{{-n, n+2} \atop \frac{3}{2}}\right|\frac{1-x}{2}\right),\;n=0,1,2,\ldots.$$
These polynomials can also be defined by their generating function
\begin{equation}
\label{genCheb2}
\frac{1}{1-2xt+t^2}=\sum_{n=0}^{\infty}U_n(x)t^n.
\end{equation}
Finally the classical Legendre (or spherical) polynomials $\left\{P_n(x)\right\}_{n=0}^{\infty}$
form another special case of the classical Jacobi polynomials. In
fact we have
$$P_n(x)=P_n^{(0,0)}(x)=\sum_{k=0}^n\frac{(n+k)!}{(n-k)!\,(k!)^2}
\left(\frac{x-1}{2}\right)^k,\;n=0,1,2,\ldots.$$
These Legendre polynomials can also be defined by their generating
function
\begin{equation}
\label{genLeg}
\frac{1}{\sqrt{1-2xt+t^2}}=\sum_{n=0}^{\infty}P_n(x)t^n.
\end{equation}
Note that the Legendre polynomials also form a special case of the
ultraspherical polynomials, since we have
\begin{equation}
\label{rel2}
P_n(x)=G_n^{(\frac{1}{2})}(x),\;n=0,1,2,\ldots.
\end{equation}

The classical Laguerre polynomials $\left\{L_n^{(\alpha)}(x)\right\}_{n=0}^{\infty}$
can be defined for all $\alpha$ and $n\in\{0,1,2,\ldots\}$ as
$$L_n^{(\alpha)}(x)=\sum_{k=0}^n(-1)^k\left({n+\alpha \atop n-k}\right)
\frac{x^k}{k!}=\sum_{k=0}^n(-1)^k\frac{(\alpha+k+1)_{n-k}}{(n-k)!}\,
\frac{x^k}{k!}.$$
The generating function for the classical Laguerre polynomials is
given by
\begin{equation}
\label{genLag}
(1-t)^{-\alpha-1}\exp\left(\frac{xt}{t-1}\right)=
\sum_{n=0}^{\infty}L_n^{(\alpha)}(x)t^n.
\end{equation}
Further we have for $n=0,1,2,\ldots$
\begin{equation}
\label{diffLag}
D^iL_n^{(\alpha)}(x)=(-1)^iL_{n-i}^{(\alpha+i)}(x),\;i=0,1,2,\ldots,n.
\end{equation}

Another family of continuous orthogonal polynomials is the one
named after Hermite. The classical Hermite polynomials
$\left\{H_n(x)\right\}_{n=0}^{\infty}$ can be defined by their
generating function
\begin{equation}
\label{genHer}
\exp\left(xt-\frac{1}{4}t^2\right)=\sum_{n=0}^{\infty}H_n(x)t^n.
\end{equation}
Here we used another normalization than in \cite{Askey}. This one
turns out to be more convenient in this paper. These classical
Hermite polynomials satisfy the orthogonality relation
$$\frac{1}{\sqrt{\pi}}\int_{-\infty}^{\infty}e^{-x^2}H_m(x)H_n(x)dx=
\frac{\delta_{mn}}{2^n\,n!},\;m,n=0,1,2,\ldots.$$
Further we have for $n=0,1,2,\ldots$
\begin{equation}
\label{diffHer}
D^iH_n(x)=H_{n-i}(x),\;i=0,1,2,\ldots,n.
\end{equation}
From the generating function it follows that
\begin{equation}
\label{nulHer}
H_{2n+1}(0)=0\;\textrm{ and }\;
H_{2n}(0)=\frac{(-1)^n}{2^{2n}\,n!},\;n=0,1,2,\ldots.
\end{equation}
Further we will use the kernels
\begin{equation}
\label{kernHer}
K_n(x,y)=\sum_{k=0}^n2^k\,k!\,H_k(x)H_k(y),\;n=0,1,2,\ldots.
\end{equation}
By using (\ref{nulHer}) we easily find that for $n=0,1,2,\ldots$
\begin{equation}
\label{kernHer1}
K_{2n+1}(x,0)=K_{2n}(x,0)=\sum_{k=0}^n(-1)^k\frac{(2k)!}{k!}H_{2k}(x)
\end{equation}
and
\begin{equation}
\label{kernHer2}
K_{2n+1}(0,0)=K_{2n}(0,0)=\sum_{k=0}^n\frac{(2k)!}{2^{2k}\,(k!)^2}
=\sum_{k=0}^n\frac{(\frac{1}{2})_k}{k!}=\frac{(\frac{3}{2})_n}{n!}.
\end{equation}

Finally we will consider the discrete orthogonal polynomials named
after Meixner and Charlier. We choose normalizations different
from those in \cite{Askey}. The classical Meixner polynomials
$\left\{M_n^{(\beta)}(x;c)\right\}_{n=0}^{\infty}$ can be defined
by their generating function
\begin{equation}
\label{genMeixner}
\left(1-\frac{t}{c}\right)^x\left(1-t\right)^{-x-\beta}=
\sum_{n=0}^{\infty}M_n^{(\beta)}(x;c)t^n.
\end{equation}
The Meixner polynomials are connected to the classical Jacobi
polynomials in the following way
\begin{equation}
\label{rel3}
M_n^{(\beta)}(x;c)=P_n^{(\beta-1,-n-\beta-x)}\left(\frac{2-c}{c}\right),
\;n=0,1,2,\ldots.
\end{equation}
The classical Charlier polynomials $\left\{C_n^{(a)}(x)\right\}_{n=0}^{\infty}$
can also be defined by their generating function
\begin{equation}
\label{genChar}
e^{-at}\left(1+t\right)^x=\sum_{n=0}^{\infty}C_n^{(a)}(x)t^n.
\end{equation}

\section{Some inversion formulas}

In \cite{Charlier} we observed that the generating function
(\ref{genChar}) implies that
$$1=e^{-at}(1+t)^xe^{at}(1+t)^{-x}=\sum_{n=0}^{\infty}
\left(\sum_{k=0}^nC_k^{(a)}(x)C_{n-k}^{(-a)}(-x)\right)t^n,$$
which implies that
$$\sum_{k=0}^nC_k^{(a)}(x)C_{n-k}^{(-a)}(-x)=
\left\{\begin{array}{ll}
1, & n=0\\
0, & n=1,2,3\ldots
\end{array}\right.$$
or
\begin{equation}
\label{invChar}
\sum_{k=j}^iC_{i-k}^{(-a)}(-x)C_{k-j}^{(a)}(x)=
\delta_{ij},\;j\le i,\;i,j=0,1,2,\ldots.
\end{equation}
As already indicated in \cite{Invjac} this formula (\ref{invChar}) can
be interpreted as follows. If we define the matrix $T=(t_{ij})_{i,j=0}^n$
with entries
$$t_{ij}=\left\{\begin{array}{ll}
C_{i-j}^{(a)}(x), & j\le i\\
0, & j>i,
\end{array}\right.$$
then this matrix $T$ is a triangular matrix with determinant $1$
and the inverse $U$ of $T$ is given by $T^{-1}=U=(u_{ij})_{i,j=0}^n$
with entries
$$u_{ij}=\left\{\begin{array}{ll}
C_{i-j}^{(-a)}(-x), & j\le i\\
0, & j>i.
\end{array}\right.$$
Therefore we call (\ref{invChar}) an inversion formula.

In the same way we find from the generating function
(\ref{genLag}) for the classical Laguerre polynomials
$$\sum_{k=j}^iL_{i-k}^{(-\alpha-2)}(-x)L_{k-j}^{(\alpha)}(x)=
\delta_{ij},\;j\le i,\;i,j=0,1,2,\ldots.$$
However, in view of (\ref{diffLag}) this inversion formula cannot
be used to solve the systems of the equations of the form
(\ref{sysLag}). In \cite{Bavinck} H.~Bavinck observed that it also
follows from the generating function (\ref{genLag}) that
\begin{eqnarray*}
(1-t)^{i-j-1}&=&(1-t)^{-\alpha-j-1}\exp\left(\frac{xt}{t-1}\right)
(1-t)^{\alpha+i}\exp\left(\frac{-xt}{t-1}\right)\\
&=&\sum_{n=0}^{\infty}\left(\sum_{k=0}^nL_k^{(\alpha+j)}(x)
L_{n-k}^{(-\alpha-i-1)}(-x)\right)t^n
\end{eqnarray*}
which implies, by comparing the coefficients of $t^{i-j}$ on both sides, that
$$\sum_{k=0}^{i-j}L_k^{(\alpha+j)}(x)L_{i-j-k}^{(-\alpha-i-1)}(-x)=\delta_{ij},
\;j\le i,\;i,j=0,1,2,\ldots,$$
which is equivalent to (\ref{invLag}). This inversion formula
implies that the system of equations (\ref{sysLag}) has the
unique solution given by (\ref{oplLag}).

\section{More (inversion) formulas}

Applying the method described in the preceding section to the
generating function (\ref{genCheb2}) for the Chebyshev
polynomials of the second kind and the generating function
(\ref{genLeg}) for the Legendre polynomials we obtain
\begin{eqnarray*}
\sum_{n=0}^{\infty}U_n(x)t^n&=&\frac{1}{1-2xt+t^2}=
\frac{1}{\sqrt{1-2xt+t^2}}\frac{1}{\sqrt{1-2xt+t^2}}\\
&=&\sum_{k=0}^{\infty}P_k(x)t^k\sum_{m=0}^{\infty}P_m(x)t^m=
\sum_{n=0}^{\infty}\left(\sum_{k=0}^nP_k(x)P_{n-k}(x)\right)t^n,
\end{eqnarray*}
which implies that
$$\sum_{k=0}^nP_k(x)P_{n-k}(x)=U_n(x),\;n=0,1,2,\ldots.$$
Another interesting formula of this kind can be found by using the
generating function (\ref{genCheb2}) for the Chebyshev polynomials
of the second kind and the generating function (\ref{genCheb1})
for the Chebyshev polynomials of the first kind. In fact, we have
\begin{eqnarray*}
\sum_{n=0}^{\infty}U_n(x)t^n&=&\frac{1}{1-2xt+t^2}=
\frac{1}{1-xt}\frac{1-xt}{\sqrt{1-2xt+t^2}}\\
&=&\sum_{k=0}^{\infty}x^kt^k\sum_{m=0}^{\infty}T_m(x)t^m=
\sum_{n=0}^{\infty}\left(\sum_{k=0}^nx^kT_{n-k}(x)\right)t^n,
\end{eqnarray*}
which implies that
\begin{equation}
\label{relCheb}
\sum_{k=0}^nx^kT_{n-k}(x)=U_n(x),\;n=0,1,2,\ldots.
\end{equation}

As before we can use the generating function (\ref{genultra}) for
the ultraspherical polynomials to obtain the inversion formula
\begin{equation}
\label{invultra}
\sum_{k=j}^iG_{i-k}^{(-\lambda)}(x)G_{k-j}^{(\lambda)}(x)=
\delta_{ij},\;j\le i,\;i,j=0,1,2,\ldots.
\end{equation}

In view of (\ref{rel2}) the special (limit) case $\lambda=\frac{1}{2}$
should lead to an inversion formula for the Legendre polynomials.
If we define for every positive integer $N$ the matrix
$A=(a_{ij})_{i,j=1}^N$ with entries
$$a_{ij}=\left\{\begin{array}{ll}
P_{i-j}(x), & j\le i\\
0, & j>i,
\end{array}\right.$$
then this matrix is a triangular matrix with determinant $1$ and
hence invertible. Now we have $G_0^{(\lambda)}(x)=1$,
$G_1^{(\lambda)}(x)=2\lambda x\;\rightarrow\;-x$ for
$\lambda\rightarrow -\frac{1}{2}$ and for $n=2,3,4,\ldots$
$$G_n^{(\lambda)}(x)=\frac{(2\lambda)_n}{(\lambda+\frac{1}{2})_n}
P_n^{(\lambda-\frac{1}{2},\lambda-\frac{1}{2})}(x)
\;\rightarrow\;\frac{-2}{n-1}P_n^{(-1,-1)}(x)\;\textrm{ for }\;
\lambda\rightarrow -\frac{1}{2}.$$
Now we have by using (\ref{defJac}) for $n=2,3,4,\ldots$
\begin{equation}
B_n(x):=\frac{-2}{n-1}P_n^{(-1,-1)}(x)=
\frac{1}{n}(1-x)P_{n-1}^{(1,-1)}(x).
\end{equation}
Hence, the inverse $A^{-1}=B=(b_{ij})_{i,j=1}^N$ is given by
$$b_{ij}=\left\{\begin{array}{ll}
0, & i<j\\
1, & i=j\\
-x & i=j+1\\
B_{i-j}(x), & i\ge j+2.
\end{array}\right.$$

In case of the Chebyshev polynomials of the second kind we can
obtain an inversion formula as follows. If we define for every
positive integer $N$ the matrix $A=(a_{ij})_{i,j=1}^N$ with entries
$$a_{ij}=\left\{\begin{array}{ll}
U_{i-j}(x), & j\le i\\
0, & j>i,
\end{array}\right.$$
then this matrix is a triangular matrix with determinant $1$ and
hence invertible. It is not difficult to show that its inverse
$A^{-1}=B=(b_{ij})_{i,j=1}^N$ is given by
$$b_{ij}=\left\{\begin{array}{ll}
1, & i=j\\
-2x, & i=j+1\\
1, & i=j+2\\
0, & \textrm{otherwise.}
\end{array}\right.$$
This can be shown by writing
$$BA=C=(c_{ij})_{i,j=1}^N\;\textrm{ with }\;c_{ij}=\sum_{k=1}^Nb_{ik}a_{kj}$$
and showing that $C=I$, the identity matrix. This is done by using
the well-known relation
$$U_n(x)-2xU_{n+1}(x)+U_{n+2}(x)=0,\;n=0,1,2,\ldots.$$

In case of the Chebyshev polynomials of the first kind we consider
the matrix $A=(a_{ij})_{i,j=1}^N$ for every positive integer $N$
with entries
$$a_{ij}=\left\{\begin{array}{ll}
T_{i-j}(x), & j\le i\\
0, & j>i.
\end{array}\right.$$
Then this matrix is also a triangular matrix with determinant $1$
and hence invertible. The inverse $A^{-1}=B=(b_{ij})_{i,j=1}^N$
is given by
$$b_{ij}=\left\{\begin{array}{ll}
0, & i<j\\
1, & i=j\\
-x, & i=j+1\\
x^{i-j-2}(1-x^2), & i\ge j+2.
\end{array}\right.$$
This can also be shown by writing
$$BA=C=(c_{ij})_{i,j=1}^N\;\textrm{ with }\;c_{ij}=\sum_{k=1}^Nb_{ik}a_{kj}$$
and showing that $C=I$, the identity matrix. This is done by using
the formula (\ref{relCheb}) and the well-known relation
$$(1-x^2)U_n(x)-xT_{n+1}(x)+T_{n+2}(x)=0,\;n=0,1,2,\ldots.$$

The generating function (\ref{genMP}) can also be used to find
inversion formulas involving Meixner-Pollaczek polynomials. In
fact we have
$$\sum_{k=j}^iP_{i-k}^{(-\lambda)}(-x;\phi)P_{k-j}^{(\lambda)}(x;\phi)
=\delta_{ij},\;j\le i,\;i,j=0,1,2,\ldots$$
or
$$\sum_{k=j}^iP_{i-k}^{(-\lambda)}(x;-\phi)P_{k-j}^{(\lambda)}(x;\phi)
=\delta_{ij},\;j\le i,\;i,j=0,1,2,\ldots.$$

By using the generating function (\ref{genMeixner}) for the
Meixner polynomials we find the inversion formula
\begin{equation}
\label{invMeixner}
\sum_{k=j}^iM_{i-k}^{(-\beta)}(-x;c)M_{k-j}^{(\beta)}(x;c)=
\delta_{ij},\;j\le i,\;i,j=0,1,2,\ldots.
\end{equation}
We remark that this inversion formula is different from the one
obtained in \cite{Meixner}. See also \cite{Meixner2} for an
application of that inversion formula.

Note that the generating function (\ref{genHer}) for the classical
Hermite polynomials implies that
$$1=\exp\left(xt-\frac{1}{4}t^2\right)\exp\left(-xt+\frac{1}{4}t^2\right)
=\sum_{n=0}^{\infty}\left(\sum_{k=0}^nH_k(x)H_{n-k}(ix)i^{n-k}\right)t^n,$$
which implies that
$$\sum_{k=0}^nH_k(x)H_{n-k}(ix)i^{n-k}=
\left\{\begin{array}{ll}
1, & n=0\\
0, & n=1,2,3\ldots.
\end{array}\right.$$
In view of (\ref{diffHer}) this formula can be used as follows.
A system of equations of the form
\begin{equation}
\label{sysHer}
F_n(x)=\sum_{k=1}^{\infty}a_k(x)D^kH_n(x),\;n=1,2,3,\ldots,
\end{equation}
where the coefficients $\left\{a_k(x)\right\}_{k=1}^{\infty}$ are
polynomials which are independent of $n$, has the unique solution
\begin{equation}
\label{oplHer}
a_k(x)=\sum_{j=1}^ki^{k-j}H_{k-j}(ix)F_j(x),\;k=1,2,3,\ldots.
\end{equation}

\section{Inversion formulas involving Jacobi polynomials}

In \cite{Invjac} we have found the inversion formula
(\ref{invJac}) involving Jacobi polynomials. As mentioned before
this formula was found in a completely different way. The
well-known generating function for the classical Jacobi
polynomials has a different structure, which implies that the
method used before cannot be used in that case. In \cite{Invjac}
we proved that for $n=0,1,2,\ldots$ we have
$$\sum_{k=0}^n\frac{\alpha+\beta+2k+1}{(\alpha+\beta+k+1)_{n+1}}
P_k^{(\alpha,\beta)}(x)P_{n-k}^{(-n-\alpha-1,-n-\beta-1)}(y)=
\frac{1}{n!}\left(\frac{x-y}{2}\right)^n.$$
Now $y=x$ leads to the inversion formula (\ref{invJac}). If $y=-x$
this leads to a formula which was used in \cite{Madrid} in the
case that $\beta=\alpha$.

By using the relation (\ref{rel3}) between the Meixner and the Jacobi
polynomials we find from the inversion formula (\ref{invMeixner})
for the Meixner polynomials that
$$\sum_{k=j}^iP_{i-k}^{(-\alpha-2,-\beta-i+k)}(x)P_{k-j}^{(\alpha,\beta-k+j)}(x)
=\delta_{ij},\;j\le i,\;i,j=0,1,2,\ldots.$$

Another inversion formula involving Jacobi polynomials can be
obtained from the inversion formula (\ref{invultra}) for the
ultraspherical polynomials. By using (\ref{rel1}) and after setting
$\lambda=\alpha+\frac{1}{2}$ this leads to
\begin{eqnarray*}
& &\sum_{k=j}^i\frac{(-2\alpha-1)_{i-k}}{(-\alpha)_{i-k}}
\frac{(2\alpha+1)_{k-j}}{(\alpha+1)_{k-j}}\\
& &{}\hspace{1cm}\times P_{i-k}^{(-\alpha-1,-\alpha-1)}(x)P_{k-j}^{(\alpha,\alpha)}(x)=
\delta_{ij},\;j\le i,\;i,j=0,1,2,\ldots.
\end{eqnarray*}

\section{Applications to differential equations}

In this section we will investigate the generalized Hermite polynomials
$\left\{H_n^M(x)\right\}_{n=0}^{\infty}$ which are orthogonal on
the real line with respect to the weight function
$$w(x)=\frac{1}{\sqrt{\pi}}e^{-x^2}+M\delta(x),\;M\ge 0.$$
In \cite{Opsym} these generalized Hermite polynomials are called
special (linear) perturbations of the classical Hermite
polynomials. They can be represented in terms of the kernels
(\ref{kernHer}) as (see \cite{Opsym})
$$H_n^M(x)=H_n(x)+MQ_n(x),\;n=0,1,2,\ldots,$$
where $Q_0(x)=0$ and
$$Q_n(x)=\left|\begin{array}{cc}
H_n(x) & K_{n-1}(x,0)\\
H_n(0) & K_{n-1}(0,0)
\end{array}\right|=\sum_{k=0}^nq_{n,k}H_k(x),\;n=1,2,3,\ldots$$
with, by using (\ref{kernHer1}), for $n=1,2,3,\ldots$
$$q_{n,n}=K_{n-1}(0,0)\;\textrm{ and }\;q_{n,k}=-2^k\,k!\,H_k(0)H_n(0),
\;k=0,1,2,\ldots,n-1.$$
In \cite{Opsym} it is shown that these generalized Hermite
polynomials satisfy a differential equation of the form
\begin{equation}
\label{dvHer}
M\sum_{k=1}^{\infty}a_k(x)y^{(k)}(x)+y''(x)-2xy'(x)+(2n+M\alpha_n)y(x)=0,
\end{equation}
where the coefficients $\left\{a_k(x)\right\}_{k=1}^{\infty}$ are
polynomials with degree$[a_k(x)]\le k,\;k=1,2,3,\ldots$ which are
independent of $n$. Moreover it is shown that the 'eigenvalue'
parameters $\left\{\alpha_{2n+1}\right\}_{n=0}^{\infty}$ can be
chosen arbitrarily,
$$\alpha_0=0\;\textrm{ and }\;\alpha_{2n}=
\sum_{j=1}^n(\lambda_{2j}-\lambda_{2j-2})q_{2j,2j},\;n=1,2,3\ldots,$$
where $\lambda_n=2n,\;n=0,1,2,\ldots$. Hence,
$\lambda_{2j}-\lambda_{2j-2}=4,\;j=1,2,3,\ldots$ and by using (\ref{kernHer2})
$$q_{2j,2j}=K_{2j-1}(0,0)=\frac{(\frac{3}{2})_{j-1}}{(j-1)!},\;j=1,2,3,\ldots.$$
This implies that
$$\alpha_{2n}=4\sum_{j=1}^n\frac{(\frac{3}{2})_{j-1}}{(j-1)!}=
4\sum_{k=0}^{n-1}\frac{(\frac{3}{2})_k}{k!}=4\frac{(\frac{5}{2})_{n-1}}{(n-1)!},
\;n=1,2,3,\ldots.$$
In order to find the coefficients $\left\{a_k(x)\right\}_{k=1}^{\infty}$
we set $y(x)=H_n^M(x)=H_n(x)+MQ_n(x)$ in the differential equation
(\ref{dvHer}) and view the left-hand side as a polynomial in $M$.
Then the coefficients of this polynomial must vanish, hence
$$\sum_{k=1}^{\infty}a_k(x)D^kH_n(x)=-\alpha_nH_n(x)-Q_n''(x)+2xQ_n'(x)-2nQ_n(x),
\;n=0,1,2,\ldots$$
and
$$\sum_{k=1}^{\infty}a_k(x)D^kQ_n(x)=-\alpha_nQ_n(x),\;n=0,1,2,\ldots.$$
Since we have, by using (\ref{nulHer}) and (\ref{kernHer2}),
$$Q_{2n+1}(x)=K_{2n}(0,0)H_{2n+1}(x)=\frac{(\frac{3}{2})_n}{n!}H_{2n+1}(x),
\;n=0,1,2,\ldots$$
and
$$\frac{(\frac{3}{2})_n}{n!}\ne 0,\;n=0,1,2,\ldots$$
both systems of equations lead to
\begin{equation}
\label{Her1}
\sum_{k=1}^{2n+1}a_k(x)H_{2n+1-k}(x)=-\alpha_{2n+1}H_{2n+1}(x),\;n=0,1,2,\ldots.
\end{equation}
Further we have
$$Q_{2n}(x)=\sum_{k=0}^nq_{2n,2k}H_{2k}(x),\;n=1,2,3,\ldots,$$
where, by using (\ref{kernHer2}),
$$q_{2n,2n}=K_{2n-1}(0,0)=\frac{(\frac{3}{2})_{n-1}}{(n-1)!},\;n=1,2,3,\ldots$$
and, by using (\ref{nulHer}), for $k=0,1,2,\ldots,n-1$ and $n=1,2,3,\ldots$
$$q_{2n,2k}=-2^{2k}(2k)!\,H_{2k}(0)H_{2n}(0)
=\frac{(-1)^{n+k+1}(2k)!}{2^{2n}\,n!\,k!}.$$
Now we have for $n=1,2,3,\ldots$
\begin{eqnarray*}
Q_{2n}''(x)-2xQ_{2n}'(x)&=&
\sum_{k=0}^nq_{2n,2k}\left[H_{2k}''(x)-2xH_{2k}'(x)\right]\\
&=&\sum_{k=0}^nq_{2n,2k}\left[-4kH_{2k}(x)\right].
\end{eqnarray*}
Hence, for $n=1,2,3,\ldots$ we obtain
\begin{eqnarray*}
-\alpha_{2n}H_{2n}(x)-Q_{2n}''(x)+2xQ_{2n}'(x)-4nQ_{2n}(x)\\
=-\alpha_{2n}H_{2n}(x)-4\sum_{k=0}^n(n-k)q_{2n,2k}H_{2k}(x),
\end{eqnarray*}
which leads for $n=1,2,3,\ldots$ to
\begin{equation}
\label{Her2}
\sum_{k=1}^{2n}a_k(x)H_{2n-k}(x)=
-\alpha_{2n}H_{2n}(x)-4\sum_{k=0}^n(n-k)q_{2n,2k}H_{2k}(x).
\end{equation}
Hence, with (\ref{Her1}) and (\ref{Her2}) we have found that
$$\sum_{k=1}^{\infty}a_k(x)D^kH_n(x)=\sum_{k=1}^na_k(x)H_{n-k}(x)=F_n(x),\;n=1,2,3,\ldots,$$
where
$$\left\{\begin{array}{l}
F_{2n+1}(x)=-\alpha_{2n+1}H_{2n+1}(x),\;n=0,1,2,\ldots\\
F_{2n}(x)=-\alpha_{2n}H_{2n}(x)-4\displaystyle\sum_{k=0}^n(n-k)q_{2n,2k}H_{2k}(x),\;n=1,2,3,\ldots.
\end{array}\right.$$
This system of equations has the form (\ref{sysHer}). So we can
use (\ref{oplHer}) to conclude that
$$a_k(x)=\sum_{j=1}^ki^{k-j}H_{k-j}(ix)F_j(x),\;k=1,2,3,\ldots.$$

We emphasize that these generalized Hermite polynomials
are orthogonal with respect to a weight function consisting of the
classical Hermite weight function and a Dirac delta distribution
at the origin. Therefore these generalized Hermite polynomials
could be considered as Krall-Hermite polynomials, but these are
quite different from the Krall-Hermite polynomials considered in
\cite{Haine} which are not orthogonal.

Finally, in \cite{Kwon} it is shown that these generalized Hermite
polynomials cannot satisfy a finite order differential equation of
the form (\ref{dvHer}).

\end{document}